\documentclass[11pt]{article}

\usepackage{amssymb,latexsym,amsmath}

\usepackage{graphicx}

\begin{document}

\newcommand{\bfi}{\bfseries\itshape}

\makeatletter

\@addtoreset{figure}{section}

\def\thefigure{\thesection.\@arabic\c@figure}

\def\fps@figure{h, t}

\@addtoreset{table}{bsection}

\def\thetable{\thesection.\@arabic\c@table}

\def\fps@table{h, t}

\@addtoreset{equation}{section}

\def\theequation{\thesubsection.\arabic{equation}}

\makeatother

\newtheorem{thm}{Theorem}[section]

\newtheorem{prop}[thm]{Proposition}

\newtheorem{lema}[thm]{Lemma}

\newtheorem{cor}[thm]{Corollary}

\newtheorem{defi}[thm]{Definition}

\newtheorem{rk}[thm]{Remark}

\newtheorem{exempl}{Example}[section]

\newenvironment{exemplu}{\begin{exempl}  \em}{\hfill $\surd$

\end{exempl}}

\newcommand{\comment}[1]{\par\noindent{\raggedright\texttt{#1}

\par\marginpar{\textsc{Comment}}}}

\newcommand{\todo}[1]{\vspace{5 mm}\par \noindent \marginpar{\textsc{ToDo}}\framebox{\begin{minipage}[c]{0.95 \textwidth}

\tt #1 \end{minipage}}\vspace{5 mm}\par}

\newcommand{\ea}{\mbox{{\bf a}}}

\newcommand{\eu}{\mbox{{\bf u}}}

\newcommand{\ueu}{\underline{\eu}}

\newcommand{\ueo}{\overline{u}}

\newcommand{\oeu}{\overline{\eu}}

\newcommand{\ew}{\mbox{{\bf w}}}

\newcommand{\ef}{\mbox{{\bf f}}}

\newcommand{\eF}{\mbox{{\bf F}}}

\newcommand{\eC}{\mbox{{\bf C}}}

\newcommand{\en}{\mbox{{\bf n}}}

\newcommand{\eT}{\mbox{{\bf T}}}

\newcommand{\eL}{\mbox{{\bf L}}}

\newcommand{\eR}{\mbox{{\bf R}}}

\newcommand{\eV}{\mbox{{\bf V}}}

\newcommand{\eU}{\mbox{{\bf U}}}

\newcommand{\ev}{\mbox{{\bf v}}}

\newcommand{\eve}{\mbox{{\bf e}}}

\newcommand{\uev}{\underline{\ev}}

\newcommand{\eY}{\mbox{{\bf Y}}}

\newcommand{\eK}{\mbox{{\bf K}}}

\newcommand{\eP}{\mbox{{\bf P}}}

\newcommand{\eS}{\mbox{{\bf S}}}

\newcommand{\eJ}{\mbox{{\bf J}}}

\newcommand{\eB}{\mbox{{\bf B}}}

\newcommand{\eH}{\mbox{{\bf H}}}

\newcommand{\leb}{\mathcal{ L}^{n}}

\newcommand{\eI}{\mathcal{ I}}

\newcommand{\eE}{\mathcal{ E}}

\newcommand{\hen}{\mathcal{H}^{n-1}}

\newcommand{\eBV}{\mbox{{\bf BV}}}

\newcommand{\eA}{\mbox{{\bf A}}}

\newcommand{\eSBV}{\mbox{{\bf SBV}}}

\newcommand{\eBD}{\mbox{{\bf BD}}}

\newcommand{\eSBD}{\mbox{{\bf SBD}}}

\newcommand{\ecs}{\mbox{{\bf X}}}

\newcommand{\eg}{\mbox{{\bf g}}}

\newcommand{\paromega}{\partial \Omega}

\newcommand{\gau}{\Gamma_{u}}

\newcommand{\gaf}{\Gamma_{f}}

\newcommand{\sig}{{\bf \sigma}}

\newcommand{\gac}{\Gamma_{\mbox{{\bf c}}}}

\newcommand{\deu}{\dot{\eu}}

\newcommand{\dueu}{\underline{\deu}}

\newcommand{\dev}{\dot{\ev}}

\newcommand{\duev}{\underline{\dev}}

\newcommand{\weak}{\stackrel{w}{\approx}}

\newcommand{\mild}{\stackrel{m}{\approx}}

\newcommand{\strong}{\stackrel{s}{\approx}}

\newcommand{\weakdown}{\rightharpoondown}

\newcommand{\opg}{\stackrel{\mathfrak{g}}{\cdot}}

\newcommand{\opunu}{\stackrel{1}{\cdot}}
\newcommand{\opdoi}{\stackrel{2}{\cdot}}

\newcommand{\opn}{\stackrel{\mathfrak{n}}{\cdot}}
\newcommand{\opx}{\stackrel{x}{\cdot}}

\newcommand{\tr}{\ \mbox{tr}}

\newcommand{\Ad}{\ \mbox{Ad}}

\newcommand{\ad}{\ \mbox{ad}}

\renewcommand{\contentsname}{ }

\title{Uniform refinements, topological derivative
 and a differentiation theorem in metric spaces}

\author{Marius Buliga \\
 \\
Institute of Mathematics, Romanian Academy \\
P.O. BOX 1-764, RO 70700\\
Bucure\c sti, Romania\\
{\footnotesize Marius.Buliga@imar.ro}}

\date{This version: 23.11.2009}

\maketitle


For the importance of differentiation theorems in metric spaces (starting with 
Pansu \cite{pansu} 
Rademacher type theorem in Carnot groups)  and relations
with rigidity of embeddings see the section 
1.2 in Cheeger and Kleiner paper \cite{cheegerkleiner} and its bibliographic
references.  

Here we propose another type of differentiation theorem, which does not involve 
measures. It is therefore different from Rademacher type theorems. Instead, this
differentiation theorem (and the concept of uniformly topological derivable
function) is formulated in terms of filters of the diagonal $\Delta(X)$ in the 
product space $\displaystyle X^{2}$ of a topological space $X$.

The main result of this paper is theorem \ref{trad2}, stated in short form as: 

\begin{thm}
Let $(X, d)$ and $(X', d')$ be $Lip_{2}$ metric spaces  and 
$f: X \rightarrow X'$ a surjective bi-Lipschitz function. 
Then $f$ is uniformly topologically derivable. 
\end{thm}

The constructions presented here are a continuation of the work  
\cite{buliga0}. For another direction of development of these ideas see the 
series of papers on dilatation structures  \cite{buligadil1} \cite{buligadil2}, 
\cite{buligairq} \cite{buligagr}, which study metric spaces with dilations,  
deformations of idempotent right quasigroups or normed groupoids.

\section{Training section in 9 steps}

The purpose of this section is to explain the constructions and ideas of this article in a familiar setting. 
In fact we shall use the most trivial scene to tell our story. 

\vspace{.5cm}

{\bf 1.} We are in $\mathbb{R}^{n}$, with $n \geq 2$. We shall use the Euclidean distance d given by the norm 
$$\| x \| \ = \ \left( x_{1}^{2} + ... + x_{n}^{2}\right)^{\frac{1}{2}} \  . $$
The distance between $x, y \in \mathbb{R}^{n}$ is then $d(x,y) = \| y-x\|$. 

The space $(\mathbb{R}^{n}, d)$ is a complete, locally compact, length metric space. 
Consider the dilatation based at $x$ and with coefficient 
$\varepsilon > 0$. This dilatation is the function: 
$$\delta_{\varepsilon}^{x}(y) \ = \ x + \varepsilon(y-x)$$
Notice that the geodesic joining $x$ and $y$ is just 
$$c(s) \ = \ \delta_{s}^{x}(y) \ , \ \forall s\in [0,1]$$
This is clearly a coincidence, specific to the example at hand.

\vspace{.5cm}

{\bf 2.} The topology on $\mathbb{R}^{n}$ is the collection $\tau$ of open sets in $\mathbb{R}^{n}$, with 
respect to the distance $d$. The empty set belongs to $\tau$. A non empty set $D \subset \mathbb{R}^{n}$ is open if for any $x \in D$ there is a $r> 0$ such that $B(x,r) \subset D$, where 
the set 
$$B(x,r) \ = \ \left\{ y \in \mathbb{R}^{n} \ \mid \ d(x,y)<r \right\}$$
is the (open) ball centered at $x$ and with radius $r$. 

Otherwise said, the topology $\tau$ is generated by the open balls. In our particular case the space is separated, meaning that it is enough to use a countable class of balls to generate all the topology. 

The space is separated in another, more rudimentary sense, explained further.

\vspace{.5cm}

{\bf 3.} For any $x \in \mathbb{R}^{n}$ let us define the filter of open neighbourhoods of $x$ as 
$$\mathcal{V}(x) \ = \ \left\{ D \in \tau \ \mid \ x \in D\right\}$$
 In general, definition \ref{dfilta}, a topological filter in the topology $\tau$ is a collection 
$\mathcal{V} \subset \tau$ with the properties: 
\begin{enumerate}
\item[(a)] if $A\in \mathcal{V}$ and $B \in \tau$ such that $A \subset B$ then $B \in \mathcal{V}$, 
\item[(b)] if $A\in \mathcal{V}$ and $B \in \mathcal{V}$ then $A\cap B \in \mathcal{V}$. 
\end{enumerate}

Obviously for any $x \in \mathbb{R}^{n}$ the collection $\mathcal{V}(x)$ is a topological filter. But there are 
many other filters in the topology $\tau$. 

The space $\mathbb{R}^{n}$ is separated in a topological sense ($T_{0}$ separated). This means that 
$$\forall x, y \in \mathbb{R}^{n} \ \ \ \mathcal{V}(x) \ = \ \mathcal{V}(y) \ \ \mbox{ implies } x=y$$

\vspace{.5cm}

{\bf 4.} A more convenient way to see a topological filter $\mathcal{V}$ is to look at the associated function $$\mu : \tau \rightarrow \left\{ 0,1\right\} \ \ , \ \mu(D) \ = \ \left\{ \begin{array}{cl}
1 & \mbox{, if } D \in \mathcal{V} \\
0 & \mbox{ otherwise}
\end{array} \right.$$

The axioms (a) and (b) from the point {\bf 3.}, which define the notion of
topological filter, can be restated as (see proposition \ref{prop23} ): 
\begin{enumerate}
\item[(a)] if $A, B \in \tau$ such that $A \subset B$ then $\mu(A) \leq \mu(B)$, 
\item[(b)] for any $A,B \in \tau$ we have: 
$$\mu(A\cup B) + \mu(A\cap B) \ \geq \ \mu(A) + \mu(B)$$
\end{enumerate}
We add to these points the requirement 
\begin{enumerate}
\item[(c)] $\displaystyle \mu(\mathbb{R}^{n}) \ = \ 1$, 
\end{enumerate}
because we want to exclude the trivial case $\mu (D) = 0$ for any open set $D$. 

From now on we shall use the name "filter" both for $\mathcal{V}$ and its associated function $\mu$. Sometimes we shall say that
$$\mathcal{V} \ = \ \left\{ D \in \tau \ \mid \ \mu(D) \ = \ 1 \right\}$$ 
is the support of $\mu$ and we shall write this as: 
$$\mathcal{V} \ = \ supp \ \mu \ .$$

To the filter $\mathcal{V}(x)$ we associate the function 
$$\circ(x) \ = \ \left\{   \begin{array}{cl}
1 & \mbox{, if } x \in D \\
0 & \mbox{ otherwise}
\end{array} \right.$$
Otherwise said the function $\circ(x)$ is the Dirac measure based at $x$. This is traditionally denoted 
by $\delta_{x}$, but we prefer to keep the $\delta$ symbol for dilatations, which will be very important in the sequel. Moreover, the notation $\circ(x)$ suggest the geometrical operation to replace $x$ by an infinitesimal "hole" in the environing space , which is exactly what we shall do three steps further. 

\vspace{.5cm}

{\bf 5.} Denote by $\mathcal{A}(\tau)$ the set of all filters. Endow this set with the topology of uniform convergence. Indeed, it is sufficient to take the topology induced by the inclusion 
$$\mathcal{A}(\tau) \subset [0,1]^{\tau} \ . $$ 
A non empty  open set in $V \subset \mathcal{A}(\tau)$ is characterized by the property: 
\begin{enumerate}
\item[-]for any $\mu \in V$ there exists $D \in \tau$ such that for any $\nu \in \mathcal{A}$, if $\nu(D) = 1$ then 
$\nu \in V$.  
\end{enumerate}

We leave to the reader to check that the function 
$$\circ: \mathbb{R}^{n} \rightarrow \mathcal{A}(\tau)$$
which associate to any $x$ its filter of neighbourhoods, is a continuous  and open  injection. Therefore 
we can see the space $\mathbb{R}^{n}$ as a subset of $\mathcal{A}(\tau)$. 

We could use larger classes of "filters", like the class $\mathcal{B}(\tau)$ of functions $\mu: \tau \rightarrow [0,1]$ which satisfy axioms (a), (b), (c) from the point 4. The generalization consists in the co-domain of $\mu$ which is taken to be $[0,1]$ instead of $\left\{0,1\right\}$. The advantage of this class is that is stable with respect to convex combinations. In fact $\mathcal{A}(\tau)$ is the set of extremal points of the family $\mathcal{B}(\tau)$. The convergence in $\mathcal{B}(\tau)$ is defined as for $\mathcal{A}(\tau)$, that is induced from the uniform convergence in $\displaystyle [0,1]^{\tau}$. 

\vspace{.5cm}

{\bf 6} Let us consider now a continuous   function $f: \mathbb{R}^{n} \rightarrow \mathbb{R}^{m}$. Denote by $\tau^{n}$, $\tau^{m}$ the topologies on $\mathbb{R}^{n}$, $\mathbb{R}^{m}$ respectively. 

The function $f$ is continuous if and only if for any $D \in \tau^{m}$ we have $\displaystyle  f^{-1}(D) \in 
\tau^{n}$. This allows us to define  the following extension of the function $f$ to the function
$$f* : \mathcal{A}(\tau^{n}) \rightarrow \mathcal{A}(\tau^{m})$$
which sends each $\mu \in \mathcal{A}(\tau^{n})$ to $f*\mu \in
\mathcal{A}(\tau^{m})$, given by: 
$$f*\mu(D) \ = \ \mu(f^{-1}(D))$$

Under the identification of  $\mathbb{R}^{n}$ with a subset of $\mathcal{A}(\tau^{n})$, the function $f*$ is  an extension of  $f$. Indeed, for any $x \in \mathbb{R}^{n}$ and for any non empty  open set $D \subset 
\mathbb{R}^{m}$   we have 
$$f*\circ(x)(D) \  = \ \circ(x)(f^{-1}(D)) \ = \ \left\{ \begin{array}{cl}
1 & \mbox{ , if } f(x) \in D \\
0 & \mbox{ , otherwise}
\end{array} \right. \ = \ \circ(f(x))$$

The extension $f*$ is continuous (see proposition \ref{prop26} ). 

\vspace{.5cm}

{\bf 7.} We introduce now the notion of refinement (see also the notion of "boundary map", definition 2.1, 
page 459, \cite{buliga0}) of the space $\displaystyle (\mathbb{R}^{n}, \tau^{n})$. In this paper the general definition is 
definition \ref{def27}. 

A function 
$$ x \in \mathbb{R}^{n} \ \ \mapsto \partial x \subset  \mathcal{A}(\tau^{n})$$
 is a refinement if for any $\displaystyle x \in \mathbb{R}^{n}$ the set $\partial x$ is non empty and for any 
 $\mu \in \partial x $ we have $\mu \geq \circ(x)$. 
 
 A topological definition of derivative now follows. The continuous function $\displaystyle f: \mathbb{R}^{n} \rightarrow \mathbb{R}^{m}$ is topologically derivable if  it commutes with the refinement, that is for any $\displaystyle x \in \mathbb{R}^{n}$ and $\mu \in \partial x$ we have 
 $f*\mu \in \partial f(x)$. 
 
 This notion is more interesting for homeomorphisms and for space dimensions greater than one 
 (we shall see further that "dimension" word may refer to metric dimension, not to  topological dimension).

We shall  give a refinement which justifies the name "derivative".   For any $\displaystyle x \in \mathbb{R}^{n}$, any $\displaystyle u \in \mathbb{R}^{n}$ with $\| u \| = 1$, and any $\varepsilon > 0$, 
$\lambda \in (0,1)$ we define the open set: 
$$V^{+} (x,u,\varepsilon, \lambda) \ = \ \left\{ y \not = x \ \mid \  \exists \sigma \in [0,\varepsilon]  \ ,  \ 
d(y, x+ \sigma u) < \lambda d(x,y) \right\}$$

Define the filter 
$\displaystyle \mu_{x,u}: \tau \rightarrow  \left\{0,1\right\}$ by: $\displaystyle \mu_{x,u}(D) = 1$ if 
there are $\varepsilon > 0$ and $\sigma \in (0,1)$ such that $\displaystyle V^{+}(x,u,\varepsilon,\sigma) 
\subset D$, and $\displaystyle \mu_{x,u}(D) = 0$ otherwise. 

The sets  $\displaystyle V^{+}(x,u,\varepsilon,\sigma)$ are given by: 
$$V^{+}(x,u,\varepsilon,\sigma) \ = \ \left\{ y \not = x \ \mid \ \|y-x - [0,\varepsilon] u \| \ < \ \sigma \|y-x\| \right\}$$
and they are defined for any $\varepsilon> 0$ and any $\sigma \in (0,1)$. 

By the notation 
$[0,\varepsilon] u$ we mean the segment 
$$[0,\varepsilon] u \ = \ \left\{ \lambda u \ \mid \ \lambda \in [0,\varepsilon] \right\}$$
The notation $\displaystyle \|y-x - [0,\varepsilon] u \|$ means the distance between the point $y-x$ and 
the segment $[0,\varepsilon]u$. 

The set $\displaystyle V^{+}(x,u,\varepsilon,\sigma)$ is therefore the set of all $y \not = x$ such that 
there is a $\lambda \in [0,\varepsilon]$ such that 
$$\| y - (x+ \lambda u)\| \ < \ \sigma \|y-x\|$$

This set has the shape of a cone with opening of angle $\alpha$ such that $\sin \alpha = \sigma$ , ended by a spherical cap with center approximatively  $x+\varepsilon u$ and radius approximatively  $\varepsilon/(1-\sigma)$ (we leave this for the reader to check). 

The filter $\displaystyle \mu_{x,u}$ is a way to see the direction $u$ starting from $x$. Indeed, filters induce a notion of convergence (of sequences in our case, because we have topology generated by a countable basis in our example). We describe this convergence in the next definition. 

\begin{defi}
Let $\mu$ be a filter in the topological space $(X,\tau)$, with countable basis of neighbourhoods. We say 
that a sequence $\displaystyle (x_{h})_{h \in \mathbb{N}}$ converges to $\mu$ if for any $D$ such that 
$\mu(D) = 1$ there is $n(D) \in \mathbb{N}$ such that for all $h \geq n(D)$ $\displaystyle x_{h} \in D$. This is equivalent to say that the sequence $\displaystyle (\circ(X_{h}))_{h \in \mathbb{N}}$ converges 
to $\mu$ in $\mathcal{A}(\tau)$. 
\end{defi}

We have the easy result: 

\begin{prop}
Consider the sequence   $\displaystyle (x_{h})_{h \in \mathbb{N}}$ such that for any $h \in \mathbb{N}$ 
$\displaystyle x_{h} \not = x$. This sequence  converges to $\displaystyle \mu_{x,u}$ if and 
only if $\displaystyle (x_{h})_{h \in \mathbb{N}}$ converges to $x$ and 
$$\lim_{h \rightarrow \infty} \frac{x_{h} - x}{\| x_{h} - x\|} \ = \ u$$
\end{prop}

The refinement we have in mind is then: 
$$\partial x \ = \ \left\{ \mu_{x,u} \ \mid \ \|u\| = 1\right\}$$

We have the theorem: 

\begin{thm}
Let $\displaystyle f: \mathbb{R}^{n} \rightarrow \mathbb{R}^{n}$ be a differentiable function such that for 
any $\displaystyle x \in \mathbb{R}^{n}$ the derivative $Df(x)$ is invertible. Then $f$ is derivable in the topological sense. 
\end{thm}

The proof is left to the reader. 

\vspace{.5cm}

{\bf 8.} We can see the refinement from the point of view of uniformities. Let us consider the uniformity 
associated to the euclidean distance in $\mathbb{R}^{n}$. This is the filter $u_{n}$ in the topological 
space $\displaystyle (\mathbb{R}^{n}, \tau^{n})^{2}$, given by: 
$$ \mbox{ if } \exists s> 0 \ B(s) \subset U \mbox{ then } u_{n}(U) \ = \ 1 \ , $$
where $\displaystyle B(s) \ = \ \left\{ (x,y) \in \mathbb{R}^{2n} \ \mid \ \|x-y\| < s \right\}$. 

For any $u \in \mathbb{R}^{n}$, $\|u \| = 1$, let us define the filter $\mu_{u}$ generated by the sets: 
$$V^{+}(u, \varepsilon, \mu) \ = \ \left\{ (x,y) \ \mid \ d(y, x+ [0,\varepsilon] u) < \mu d(x,y) \right\}$$
for all $\varepsilon> 0$ and $\mu \in (0,1)$. 

It can be easily checked that: 
\begin{enumerate}
\item[(a)] for any $u \in \mathbb{R}^{n}$, $\|u \| = 1$ we have $\mu_{u} \geq u_{n}$. 
\item[(b)] for any $u \in \mathbb{R}^{n}$, $\|u \| = 1$, for any $U$ such that $\mu_{u}(U) = 1$ exists 
$D$ such that $D \circ D \subset U$, where in general 
$$U \circ V \ =  \ \left\{ (x,z) \ \mid \ \exists y \ \ (x,y) \in U \ , \ (y,z) \in V \right\} \ .$$
\item[(c)] for any $u \in \mathbb{R}^{n}$, $\|u \| = 1$ we have 
$$\sigma  * \mu_{u} \ = \ \mu_{-u} $$
where $\sigma(x,y) \ = \ (y,x)$. 
\end{enumerate}

This is the definition of an uniform refinement (compare with definition
\ref{def29}). To any function $f: \mathbb{R}^{n} \rightarrow \mathbb{R}^{m}$ 
we associate the function: 
$$f^{2}: \mathbb{R}^{2n} \rightarrow \mathbb{R}^{2m} \ , \ \ f(x,y) \ = \ (f(x,x), f(x,y))$$
If the function $\displaystyle f^{2}$ transports one uniform refinement to another then we say that it 
is uniformly topologically derivable. 

An important class of uniform functions are the Lipschitz functions. We shall 
see that we can construct uniform refinements such that any bi-Lipschitz is 
uniformly derivable. 

\vspace{.5cm}

{\bf 9.} In the frame of uniform refinements we have a natural operation. Let $\mu, \nu$ filters in 
$\mathbb{R}^{2n}$ such that 
$$\circ(\Delta(\mathbb{R}^{n})) \leq \mu \ , \  \circ(\Delta(\mathbb{R}^{n})) \leq \nu \ .$$ 
Here we used the notation: 
$$\Delta(\mathbb{R}^{n}) \ = \ \left\{ (x,x) \ \mid \ x \in \mathbb{R}^{n} \right\} \ .$$

The composition $\mu \circ \nu$ is defined by: 
$$\mbox{ if } \exists A,B \ \mu(A) = 1 \ , \ \nu(B) =  1 \ , \ A \circ B \subset U \ ,  \mbox{ then } \mu\circ\nu(U) = 1 \ .$$
The operation is a priori non commutative. By chance we have:  for any $u,v \in \mathbb{R}^{n}$, 
$\|u\| = \|v\| = 1$ 
$$\mu_{u} \circ \mu_{v} \ = \ \mu_{v} \circ \mu_{u} \ .$$
This comes from the fact that the operation in $\mathbb{R}^{n}$ is commutative and from the fact that 
all filters $\mu_{u}$ are   invariant with respect to translations. 

It is natural to expect that for a general uniform refinement the composition operation is not commutative. 
This leads to "rigidity" results: suppose that we have two spaces $X$ and $Y$, each with an uniform refinement, but such that the refinement on $X$ is commutative but the refinement on $Y$ is not. Then there is no uniform function from $X$ to $Y$ which transports $\mu, \nu$ in the first refinement to 
non-commuting $\mu'$, $\nu'$. Indeed, the reason is that  the transport is a morphism with respect to filter composition. 

We shall see that  some difficult  rigidity  theorems are consequences of this simple fact.

\section{Topological derivative and uniform refinements}

\subsection{Topological derivative}

The topological derivative has been introduced in \cite{buliga0}. We see it like a natural construction in the frame of topological spaces and continuous maps. 

Let $(X,\tau)$ be a topological space. $\tau$ is the collection of open sets in $X$. 

\begin{defi}
A filter in $(X,\tau)$ is a function $\mu: \tau \rightarrow \left\{0,1\right\}$ such that: 
\begin{enumerate}
\item[(a)] $\mu(X) \ = \ 1$, 
\item[(b)] for any $A,B \in \tau$, if $A \subset B$ then $\mu(A) \leq \mu(B)$,
\item[(c)] for any $A,B \in \tau$ we have $\displaystyle \mu(A\cup B) + \mu(A\cap B) \geq \mu(A) + \mu(B)$.
\end{enumerate}
We shall denote the set of filters by $\mathcal{A}(\tau)$. 
\label{dfilta}
\end{defi}

There are other two interesting classes of filters in a weaker sense. 

\begin{defi}A function $\mu: \tau \rightarrow [0,1]$ which satisfies (a), (b),
(c) of previous definition is called a filter of type B. The set of all these filters is denoted by $\mathcal{B}(\tau)$. 

Let $\mathbb{P}(\mathbb{R})$ be the class of Borel probabilities on $\mathbb{R}$. A filter of class C is a map $\mu: \tau \rightarrow \mathbb{P}(\mathbb{R})$ such that for any measurable 
set $I \subset \mathbb{R}$ we have 
$$\mu(X)(I) \not = 0$$ 
and the map $$\displaystyle \tilde{\mu}(\cdot)(I) \ = \ \frac{1}{\mu(X)(I)} \mu(\cdot)(I)$$
is a class B filter. We denote this class of filters by $\mathcal{C}(\tau)$.  
\label{dfiltbc}
\end{defi}

Let $\mu \in \mathcal{A}(\tau)$ or $\mu \in \mathcal{B}(\tau)$ be a filter or a B class filter. The support of $\mu$ is 
$$\mbox{supp}\ \mu \ = \ \left\{ D \in \tau \ \mid \ \mu(D) = 1 \right\}$$

\begin{prop}
$\mathcal{V} \ = \ \mbox{supp} \ \mu$ is a topological filter in the usual sense, that is 
$\mathcal{V} \subset \tau$ satisfies: 
\begin{enumerate}
\item[(a')] $X \in \mathcal{V}$, 
\item[(b')] if $A \in \mathcal{V}$ and $B \in \tau$ such that $A \subset B$, then $B \in \mathcal{V}$, 
\item[(c')] if $A, B \in \mathcal{V} $ then $A\cap B \in \mathcal{V}$. 
\end{enumerate}
\label{prop23}
\end{prop}

\paragraph{Proof.} Indeed, properties (a'), (b') and (c') correspond respectively to (a), (b), (c) from definition \ref{dfilta}. It is obvious why (a') comes from (a) and (b') comes from (b). 

We prove only (c'). Take $A,B \in \mathcal{V}$. Then $\mu(A) = \mu(B) = 1$.
Because $A \subset A\cup B$ we deduce that $\mu(A \cup B) \geq \mu(A) = 1$, so
$\mu(A\cup B) = 1$. From (c), definition \ref{dfilta}, we have that: 
$$\displaystyle \mu(A\cup B) + \mu(A\cap B) \geq \mu(A) + \mu(B)$$
From previous considerations, this relations becomes: 
$$1 + \mu(A\cap B) \geq 1+1 = 2$$
therefore $\mu(A\cap B) = 1$. This shows that $A\cap B \in \mathcal{V}$. 
\quad $\blacksquare$

Let us denote by $\mathcal{A}(\tau)$ the set of all filters over $(X,\tau)$. Suppose that   $(X,\tau)$ is a  
$T_{0}$ topological space, i.e. it has the following separation property: for any $x, y \in X$, $x \not = y$, there is $V \in \tau$ such that $x \in V$ and $y \not \in V$. Then we have a canonical injection of $X$ in 
$\mathcal{A}(\tau)$, given by: 
$$  x \in X \ \mapsto \circ(x) \in \mathcal{A}(\tau) \ ,  \ \ \circ(x)(A) \ = \ \left\{ \begin{array}{ll}
1 & , \ x \in A \\ 
0 & , \ x \not \in A \end{array} \right.$$

Further we shall always suppose that $(X,\tau)$ is $T_{0}$. 

\begin{defi}
On $\mathcal{A}(\tau)$ we put the following topology, called $\tau^{e}$: a non empty set $V \subset \mathcal{A}(\tau)$ is open if for any $\mu \in V$ there is $D \in \tau$ such that $\mu(D) = 1$ and for any 
$\mu' \in \mathcal{A}(\tau)$ $\mu'(D) = 1$ implies $\mu' \in V$. 

In the case of filters of class B, let $\mathcal{B}(\tau)$ be the collection of all such filters over 
$\tau$. The topology that we consider is $\tau^{B}$: a non empty set $V \in \tau^{B}$ has the property that for any $\mu \in V$  there are $D \in \tau$ and $\varepsilon \in (0,1)$ such that $\mu(D) = 1$ and for 
any $\mu' \in \mathcal{B}(\tau)$ $\mu'(D) > 1-\varepsilon$ implies $\mu' \in V$. 
\end{defi}

The class $\mathcal{B}(\tau)$ is convex (convex combination of filters is a filter), has $\mathcal{A}(\tau)$ as the set of extremal points.

 The topology $\tau^{B}$ is the topology of uniform convergence in the class of functions defined on $\tau$  with values in $\mathbb{R}$. Moreover $\tau^{e}$ is the induced topology 
on $\mathcal{A}(\tau) \subset \mathcal{B}(\tau)$ by $\tau^{B}$. 

\begin{defi}
Let $f: (X,\tau) \rightarrow (X',\tau')$ be a continuous function. Then 
$$f*: \mathcal{A}(\tau) \rightarrow \mathcal{A}(\tau')$$
 is defined by: for all $\mu \in \mathcal{A}(\tau)$ and for all $A \in \tau'$ 
$$f*\mu (A) \ = \ \mu(f^{-1}(A))$$
Such a function can  be equally defined for class B filters. 
\end{defi}

\begin{prop}
Let $f: (X,\tau) \rightarrow (X',\tau')$ be a continuous function. 
Then $f*: \mathcal{A}(\tau) \rightarrow \mathcal{A}(\tau')$ is continuous. 
\label{prop26}
\end{prop}

\paragraph{Proof.} Let $V'$ be open in $\mathcal{A}(\tau')$, non empty. For any $\mu \in \mathcal{A}(\tau)$, if 
$f*\mu \in V'$ then  there is $D' \in \tau'$ such that $f*\mu(D') = 1$ and for any  $\lambda' \in \mathcal{A}(\tau')$ 
$\lambda'(D') = 1$ implies $\lambda' \in V'$. In particular for any $\lambda \in \mathcal{A}(\tau)$, if 
$f*\lambda(D') = 1$ then $f*\lambda \in V'$. 

Let $D = f^{-1}(D') \in \tau$. By previous reasoning we have that $\mu(D) = 1$ and for any $\lambda \in 
\mathcal{A}(\tau)$, if $\lambda(D) = 1$ then $f*\lambda \in V'$. 

This shows that 
$\displaystyle \left(f*\right)^{-1}(V')$ is open in $\mathcal{A}(\tau)$. \quad $\blacksquare$

We shall use this proposition in order to construct a derivative in the topological category. The notion that we need is refinement. 

\begin{defi}
A refinement of the $T_{0}$ topological space $(X,\tau)$ is a map
$$x \in  X \ \mapsto \ \partial(x) \subset  \mathcal{A}(\tau)$$
which has the following properties: 
\begin{enumerate}
\item[(a)] for any $\mu \in \partial(x)$ and any $D \in \tau$ we have $\mu(D) \geq \circ(x)(D)$, 
\item[(b)] for any $x \in X$ and $\mu \in \partial(x)$ there is a $D \in \tau$ such that $\mu(D) \not = \circ(x)(D)$. 
\end{enumerate}
The space $\partial X$ is then equal with the reunion of all $\partial(x)$, $x \in X$.
\label{def27}
\end{defi}

An interesting problem: what are the implications of the supplementary assumption that the refinement 
is a bundle? By this we mean: the map $\pi: \partial X \rightarrow X$ which sends $\mu \in \partial (x)$ to $x$, is a bundle. The weak assumption that we can locally trivialize the bundle can turn to be very powerful. We shall touch this subject in a future paper. This is related to the local contractibility condition.

One of the main definitions of this paper is this. 

\begin{defi}
Let $(X,\partial)$, $(Y,\partial')$ be two topological spaces with refinements $\partial$ and $\partial'$ respectively. A continuous function $f: X \rightarrow Y$ is $\partial-\partial'$ (topologically) derivable if 
$$f * \partial \ = \ \partial' \circ f$$
This means that for any $x \in X$ and $\mu \in \partial(x)$ we have $$f * \mu \in \partial'(f(x)) \ .$$ 
\end{defi}

\subsection{Uniform refinements}

Let $(X,\tau)$ be a $T_{0}$ topological space. We denote by $\Delta(X)$ the 
diagonal of $\displaystyle X^{2}$, and by $\circ (\Delta(X))$ the filter 
of neighbourhoods of $\Delta(X)$. For sets $\displaystyle A,B \subset X^{2}$, let $\displaystyle A \circ B \subset X^{2}$ be the 
set $$A \circ B \ = \ \left\{ (x,z) \in X^{2} \ \mid \ \exists y \in X \ ,\ (x, y) \in A \ \mbox{ and } (y,z) \in B \right\}$$
This composition operation defines a composition for filters in $\displaystyle
(X^{2}, \tau^{2})$.  Let $\mu, \lambda$ such filters, and $\displaystyle D \in \tau^{2}$, then  
$$\mu \circ \lambda (D) = 1 \mbox{ iff } \exists E, F \in \tau^{2} \ , \ \mu(E) = 1 \ , \ \lambda(F) = 1  \ , \ \ E\circ F \subset D$$
Let $\displaystyle \sigma: X^{2} \rightarrow X^{2}$ be the map $\sigma(x,y) = (y,x)$. 

\begin{defi}
An uniformity on $(X,\tau)$ is a filter $\Omega$ on 
$\displaystyle (X^{2}, \tau^{2})$ such that: 
\begin{enumerate}
\item[(a)]   $\Omega \geq \circ(\Delta(X))$,  
\item[(b)] for any $D$ such that $\mu(D) = 1$ 
there is $E$, $\mu(E) = 1$ and $E \circ E \subset D$. 
\item[(c)] $\displaystyle \sigma * \Omega = \Omega$.    
\end{enumerate}
\end{defi}

Remark that condition (b)  implies  (but is not equivalent to): $$\Omega \circ \Omega \geq \Omega \ .$$

An uniform refinement is a replacement of the uniformity $\Omega$ with a set of finer filters $\partial \Omega$, which satisfy a number of requirements. Because such requirements will prove to be not 
easy to satisfy in general, we shall introduce first the notion of uniform pre-refinement. 

\begin{defi}
A pre-refinement of the uniform space $(X,\tau, \Omega)$ is  a set $\partial \Omega$ of filters 
in $\displaystyle \tau^{2}$ such that for any $\mu \in \partial \Omega$ we have $\mu \geq \Omega$. 
An uniform pre-refinement $\partial \Omega$ is an uniform refinement if moreover: 
\begin{enumerate}
\item[(a)] for any $\mu \in \partial \Omega$ and  for any $D$ such that $\mu(D) = 1$ 
there is $E$, $\mu(E) = 1$ and $E \circ E \subset D$. 
\item[(b)] for any $\mu \in \partial \Omega$ we have $\sigma * \mu \in \partial \Omega$. 
\end{enumerate}
\label{def29}
\end{defi}

The derivability notion associated with uniform pre-refinements is described further. 

\begin{defi}
Let $(X,\partial \Omega)$, $(Y,\partial\Omega')$ be two uniform  spaces with refinements $\partial \Omega$ and $\partial\Omega'$ respectively. An uniformly  continuous function $f: X \rightarrow Y$ is $\partial-\partial'$ (topologically) derivable if 
$$f * \partial \Omega \ = \ \partial\Omega'$$
where $f*$ is the transport by the function $\displaystyle f^{2} : X^{2}
\rightarrow Y^{2}$, $\displaystyle f^{2}(x,y) \ = \ (f(x), f(y))$. 

This means that for any $\mu \in \partial \Omega$ we have 
$\displaystyle f^{2} * \mu \in \partial \Omega'$.  
\end{defi} 

Let now $\mu \in \partial \Omega$. For any $D$ such that $\mu(D) = 1$ and for any $x \in X$ let 
$$D_{x} \ = \ \left\{ y \not = x \ \mid \ (x,y) \in D \right\}$$
and $\displaystyle \mu_{x}$ the filter on $X$ with the property that for any $D$ such that $\mu(D) = 1$ we have $\displaystyle \mu_{x}(D_{x}) = 1$. Define now 
$$\partial(x) \ = \ \left\{ \mu_{x} \ \mid \ \mu \in \partial \Omega \right\}$$

\begin{prop}
For any given uniform pre-refinement $\partial \Omega$, the map $x \in X \mapsto \partial (x) $ is a refinement. 
\end{prop}

The easy proof is left to the reader. 

Remark the formal resemblance between the requirements on the filter $\Omega$ to be an uniformity and supplementary conditions on the uniform pre-refinement $\partial \Omega$ to be an uniform refinement. The supplementary conditions almost say that any $\mu \in \partial \Omega$ is a new uniformity, finer than the original uniformity $\Omega$.

\subsection{The snowflake functor and refinements}

This example shows that for the refinements the metric dimension is more important than the topological 
dimension. 

A well known trick to produce a bi-Lipschitz map from a Lipschitz one is to associate to the Lipschitz function $f: X \rightarrow Y$ the bi-Lipschitz function $g: X \rightarrow X \times Y$, given by: 
$$\forall x  \in X \ , \ g(x) \ = \ (x,f(x))$$
Here we take on the space $X \times Y$ the distance 
$$ d((x,y), x',y')) \  = \ d_{X}(x,x') + d_{Y}(y,y')$$

We shall take the space  $X = \mathbb{R}$ with the distance $$d(x,y) = \mid x - y \mid^{\frac{1}{m}}$$
where $m \in \mathbb{N}$, $m \geq 2$. The Hausdorff dimension of the space is therefore $m \geq 2$. 
The topological dimension is one. 

We take now the space $Y = \mathbb{R}$ with the usual distance and 
 $\displaystyle X \times Y = \mathbb{R}^{2}$ with the distance:  
$$d((x_{1}, y_{1}), (x_{2}, y_{2})) \ = \ \mid x_{1} - x_{2} \mid + \mid y_{1} - y_{2}\mid^{\frac{1}{m}} \ .$$
Let $p=p(x)$ be a polynomial of degree at most $m$, such that $p(0) = 0$. Define now the filter 
$\displaystyle \mu_{x,p}$ in $\displaystyle \mathbb{R}^{2}$  to be generated by the open sets: 
$$V^{+}(x,p, \varepsilon, \lambda) \ = \ \left\{ y \not = x \ \mid \ d(y, x+ p([0,\varepsilon]) \ < \ \lambda d(x,y) \right\}$$

We have the easy  proposition: 

\begin{prop}
Let $\displaystyle p_{i}=p_{i}(x)$, $i=1,2$  be two polynomials of degree at most $m$, such that 
$\displaystyle p_{i}(0) = 0$. We have  $\displaystyle \mu_{x,p_{1}} = \mu_{x,p_{2}}$ if and only if 
$\displaystyle p_{1} = p_{2}$. 
\end{prop}

\begin{thm}
Let $\displaystyle f: \mathbb{R} \rightarrow \mathbb{R}$ be a $m$ times differentiable function such that for any $\displaystyle x \in \mathbb{R}$ the development of order $m$ of $f$  around $x$  is nontrivial. Then $f$ is derivable in the topological sense, with respect to the refinement 
$$\partial^{m}x \ = \ \left\{ \mu_{x,p} \ \mid \ p \in \mathbb{R}[X] \ , \ 1  \leq \mbox{deg } p \leq m \ , \ p(0) = 0\right\} . $$ 
\end{thm}

\subsection{Refinements in general metric spaces}

We are in a locally compact, length metric space $(X,d)$. We shall construct a refinement of the space and a uniform refinement of the same space. We begin with the refinement.

We shall say further that the locally compact length space $(X,d)$ is of first type if  for any $x \in X$ there is a bi-Lipschitz function $c$ such that $c(0) = x$.  A shorter name is $Lip_{1}$.

Let $a<0<b$, a bi-Lipschitz curve $c: [a,b] \rightarrow X$ and $x = c(0)$. To these data we associate a
filter. For any $\varepsilon > 0$ and any $\mu \in (0,1)$ let us define the open set: 
$$V^{+}(x,c,\varepsilon, \mu) \ = \ \left\{ y \in X  \ \mid \ d(y, c([0,\varepsilon])) < \mu d(x,y) \right\} \ .$$
This family of open sets generates the filter $\partial^{+} c (x)$. We can define also the filter 
$\partial^{-}c(x)$, generated by the family of sets: 
$$V^{-}(x,c,\varepsilon, \mu) \ = \ \left\{ y \in X  \ \mid \ d(y, c([-\varepsilon, 0])) < \mu d(x,y) \right\} \ .$$
Remark that if 
$$\limsup_{\varepsilon \rightarrow 0} \sup \left\{ \frac{d(c(-s), x) + d(c(t), x)}{d(c(-s), c(t))} \ \mid \ s,t \in (0,\varepsilon) \right\} < + \infty$$
then $\partial^{+} c(x) \not = \partial^{-}c(x)$. 

Also, let $\bar{c}$ be the curve $\bar{c}(t) \ = \ c(-t)$. Then $\partial^{+}\bar{c}(x) = \partial^{-}c(x)$. Finally, let us use a bi-Lipschitz  reparametrization of  the curve c,  so that we get another bi-Lipschitz curve, $c'$, such that $c'(0) = x$. Then $\partial^{+}c' (x) = \partial^{+}c(x)$. 

We leave to the reader to prove these, if necessary. 

Let $\varepsilon> 0$ and $\mu \in (0,1)$. The set $c([0,\varepsilon])$ is compact. Therefore  $y \in V^{+}(x,c,\varepsilon, \mu)$ if and only if there is 
$\lambda \in [0,\varepsilon]$ such that 
$$d(y, c(\lambda)) < \mu d(x,y) \ .$$
From the triangle inequality we deduce that 
\begin{equation}
\frac{1}{1+\mu} d(x,c(\lambda)) < d(x,y) < \frac{1}{1-\mu} d(x,c(\lambda)) \ . 
\label{bound}
\end{equation}
As a consequence we have the following theorem, with a straightforward 
proof based on unwinding
the definitions.  

\begin{thm}
For any $x \in X$ define: 
$$\partial x \ = \ \left\{ \partial^{+} c(x) \ \mid \ c \mbox{ bi-Lip. and } c(0) = x \right\}  \ .$$
If $(X,d)$ is of first type  then $\partial$ is a refinement of the topological space $X$. Moreover, any bi-Lipschitz function $f: X \rightarrow Y$, where $X,Y$ are of first type, then $f$ is topologically derivable 
everywhere. 
\label{trad1}
\end{thm}

 We shall construct now an uniform refinement. We need more assumptions on the space $(X,d)$. We shall say that the space $(X,d)$ is of second kind (or $Lip_{2}$) if there exists a one parameter flow: 
 $$\mathcal{F}: [a,b] \times X  \rightarrow X$$ 
 such that: 
 \begin{enumerate}
 \item[(a)] for all $x\in X$ we have $\mathcal{F}(0,x) = x$, 
 \item[(b)] the family of functions $\displaystyle \left\{\mathcal{F}(\cdot, x)
 \mid x \in X \right\}$ is equicontinuous with respect to $x$, 
 \item[(c)] there is a positive function $M: [a,b] \rightarrow (0,+\infty)$, continuous in $0$, such that for any 
 $t \in [a,b]$ we have: 
 $$\limsup_{\varepsilon \rightarrow  0} \sup \left\{ d(\frac{\mathcal{F}(t,x), \mathcal{F}(t,y))}{d(x,y)} \ \mid \ 
 d(x,y)< \varepsilon \right\} \ \leq M(t) \ , $$
  $$\liminf_{\varepsilon \rightarrow  0} \inf \left\{ d(\frac{\mathcal{F}(t,x), \mathcal{F}(t,y))}{d(x,y)} \ \mid \ 
 d(x,y)< \varepsilon \right\} \ \geq \frac{1}{M(t)} \ , $$
\item[(d)] $\mathcal{F}$ is a local group: whenever $s,t, s+t \in [a,b]$ we have 
$$\mathcal{F}(s+t, x) \ = \ \mathcal{F}(s,\mathcal{F}(t,x))$$
for any $x \in X$, 
\item[(e)] A supplementary condition to be introduced soon, after some notations. 
\end{enumerate}

A sufficient condition for $(X,d)$ to be $Lip_{2}$ is to admit non-trivial bi-Lipschitz flows, with the trajectories not too wild. Remark that the trajectories of a bi-Lipschitz flow have no reason to be Lipschitz 
curves. In fact this is generically not true outside the realm of riemannian spaces. 

To any flow $\mathcal{F}$ we associate a filter $\partial^{+}\mathcal{F}$, in
$X\times X$. This filter is generated by the open sets: 
$$V^{+}(\mathcal{F}, \varepsilon, \mu) \ = \ \left\{ (x,y) \in X^{2} \ \mid \ d(y, \mathcal{F}([0,\varepsilon], x) < 
\mu d(x,y) \right\} \ , $$
for any $\varepsilon > 0$ and $\mu \in (0,1)$. We can also define the filter $\partial^{-} \mathcal{F}$, generated by all open sets 
$$V^{-}(\mathcal{F}, \varepsilon, \mu) \ = \ \left\{ (x,y) \in X^{2} \ \mid \ d(y, \mathcal{F}([-\varepsilon,0], x) < 
\mu d(x,y) \right\} \ , $$
for any $\varepsilon > 0$ and $\mu \in (0,1)$.

As previously, if we define $\bar{\mathcal{F}} (s,x) \ = \ 
\mathcal{F}(-s,x)$ then we have $\partial^{+} \bar{\mathcal{F}} \ = \ \partial^{-}\mathcal{F}$. 

The supplementary condition (e) can now be formulated: 

\begin{enumerate}
\item[(e)] There is a constant $C>0$ such that for any sufficiently small $\varepsilon', \varepsilon" > 0$, 
$\mu', \mu" \in (0,1)$ and $x,y,z \in X$, 
$$(y,z) \in V^{+}(\mathcal{F}, \varepsilon', \mu') \mbox{ and } (y,x) \in V^{+}(\bar{\mathcal{F}}, \varepsilon", \mu")$$
implies that $d(x,y) + d(y,z) \ \leq \ C d(x,z)$. 
\end{enumerate}

Denote by $Lipflow(X,d)$ the collection of flows $\mathcal{F}$ which satisfy conditions (a), ... ,(e). We have the theorem: 

\begin{thm}
(a) Let $(X,d)$ be $Lip_{2}$. Then 
$$\partial X^{2} \ = \ \left\{ \partial^{+} \mathcal{F} \ \mid \ \mathcal{F} \in \ Lipflow(X,d) \right\}$$
is an uniform refinement. 

(b) Moreover, let $(X, d)$ and $(X', d')$ be metric spaces of second kind and 
$f: X \rightarrow X'$ a surjective bi-Lipschitz function. 
Then $f$ is uniformly topologically derivable with respect to the uniform
refinements constructed at point (a). 
\label{trad2}
\end{thm}

\paragraph{Proof.} 
We have to prove that $\partial X^{2}$ is a uniform refinement. 

{\bf Step 1.} We start by showing that 
$\partial^{+} \mathcal{F} \geq \circ(\Delta(X))$. For this we have to prove the following: 
for any $\varepsilon'> 0$ there are $\varepsilon > 0$  and $\mu \in (0,1)$ such that: 
$$V^{+}(\mathcal{F}, \varepsilon, \mu) \subset B(\varepsilon') \ , $$
where we use the notation: 
$$B(\varepsilon) \ = \ \left\{ (x,y) \ \mid \ d(x,y) < \varepsilon \right\}  \ .$$
We shall use an inequality similar to \eqref{bound}. Choose $\mu \in (0,1)$ arbitrary. We have still to choose $\varepsilon > 0$. For this let $\displaystyle (x,y) \in V^{+}(\mathcal{F}, \varepsilon, \mu)$. This is equivalent to the existence of $\lambda \in [0,\varepsilon]$ such that: 
$$d(y, \mathcal{F}(\lambda, x)) < \mu d(x,y)$$
This implies, by the triangle inequality: 
$$d(x,y) \leq d(x, \mathcal{F}(\lambda, x)) + d(y, \mathcal{F}(\lambda, x)) < d(x, \mathcal{F}(\lambda, x)) + \mu d(x,y)$$
therefore we have: 
$$d(x,y) < \frac{1}{1-\mu} d(x, \mathcal{F}(\lambda, x)) \ .$$
Now we use conditions (a),  (b) on $\mathcal{F}$, to get:  there is $\varepsilon_{0} > 0$ such that for 
any $0< \lambda < \varepsilon_{0}$ we have: 
$$d(x, \mathcal{F} (\lambda, x)) \ \leq \ (1-\mu) \varepsilon' \ .$$
It is therefore sufficient to choose $\displaystyle \varepsilon \ = \ \frac{\varepsilon_{0}}{2}$.  

{\bf Step 2.} We have the following lemma: 

\begin{lema}
Given $\varepsilon'> 0$, $\mu' \in (0,1)$, sufficiently small, there are $\varepsilon"> 0$, $\mu" \in 
(0,1)$ such that 
$$(x,y) \in V^{+}(\bar{\mathcal{F}}, \varepsilon", \mu") \ \mbox{ implies } (y,x) \in V^{+}(\mathcal{F}, \varepsilon', \mu')  \ .$$
\label{lemacon}
\end{lema}

We postpone the proof of the lemma. We shall use it twice, first time to see that an immediate consequence is that $\sigma * \partial^{+}\mathcal{F} \ = \ \partial^{+}\bar{\mathcal{F}}$. 

{\bf Step 3.} Let $\varepsilon > 0$ and $\mu \in (0,1)$. We want to prove that there are $\varepsilon' > 0$ 
and $\mu' \in (0,1)$ such that: 
$$V^{+}(\mathcal{F}, \varepsilon', \mu') \circ V^{+}(\mathcal{F}, \varepsilon', \mu') \ \subset \ V^{+}(\mathcal{F}, \varepsilon, \mu) \ .$$
We start with $\varepsilon'$, $\mu'$ arbitrary and we look for sufficient conditions so that the previous inclusion happens. 

A pair $(x,z)$ belongs to $\displaystyle V^{+}(\mathcal{F}, \varepsilon', \mu')
\circ V^{+}(\mathcal{F}, \varepsilon', \mu')$ if and only if there are $\lambda_{1}, \lambda_{2} \in [0,\varepsilon']$ and $y  \in X$  such that 
$$d(y , \mathcal{F}(\lambda_{1}, x) < \mu' d(y,x) \quad \mbox{ and } \quad 
d(z , \mathcal{F}(\lambda_{2}, y) < \mu' d(z,y) \ . $$
Then we have, using conditions (a), (c), (d): 
$$d(z, \mathcal{F}([0, 2\varepsilon'], x)) \ \leq \ d(z, \mathcal{F}(\lambda_{1}+\lambda_{2}, x)) \ = $$
$$= \ d(z, \mathcal{F}(\lambda_{2}, \mathcal{F}(\lambda_{1}, x))) \ \leq \ d(z, \mathcal{F}(\lambda_{2}, y)) \ + \ d(\mathcal{F}(\lambda_{2}, y) , \mathcal{F}(\lambda_{2}, \mathcal{F}(\lambda_{1}, x))) \ < $$
$$< \ \mu' d(z,y) \ + \ 2 M(\lambda_{2}) \mu' d(y,x) \ .$$
The lemma \ref{lemacon} and condition (e) tell us that there is $C> 0$ such that  for sufficiently small $\varepsilon' > 0$ and 
$\mu' \in (0,1)$  we have: 
$$\mu' d(z,y) \ + \ 2 M(\lambda_{2}) \mu' d(y,x) \  \leq \ \max \left\{ \mu' , M(\lambda_{2}) \mu' \right\} C d(x,z) .$$
Remark that $M(0) = 1$, therefore we can choose $\varepsilon'> 0$ and $mu' \in (0,1)$ so small that 
\begin{enumerate}
\item[-] $2 \varepsilon' < \varepsilon$, 
\item[-] $\displaystyle \sup \left\{ M(\lambda) \mid \lambda \in [0,\varepsilon'] \right\} \leq 2$, 
\item[-] $4 \mu' C < \mu$. 
\end{enumerate}
This finishes the step 3 of the proof.  We proved that $\displaystyle X^{2}$ is a uniform refinement. 

{\bf Step4. } We are left with the proof of the topologically uniform derivability. For this we define first an action of bi-Lipschitz maps on  $Lipflow(X,d)$. 

\begin{lema}
Let $(X,d)$, $(X',d')$ be  $Lip_{2}$ spaces and $f: X \rightarrow X'$ a surjective bi-Lipschitz map. The following function is then well defined: 
$$f* : Lipflow(X,d) \rightarrow Lipflow(X',d') \ , \ \ f*\mathcal{F}(t,x) \ = \
f(\mathcal{F}(t, f^{-1}(x))) \ .$$
Moreover, if $(X,d)$ is $Lip_{2}$,  $(X',d')$ is locally compact length metric space and it exists a 
surjective  bi-Lipschitz map $f: X \rightarrow X'$, then $(X',d')$ is $Lip_{2}$. 
\label{lem2}
\end{lema}

Next we describe how $f$ acts on filters. Recall that for filters in $X^{2}$ we defined the action of 
$f: X \rightarrow Y$, uniformly continuous, as the action of $f^{2}: X \times X \rightarrow Y \times Y$, 
$f^{2}(x,y) \ = \ (f(x), f(y))$. 

\begin{lema}
Let $(X,d), (X',d')$ Be $Lip_{2}$ spaces and $f: X \rightarrow X'$ bi-Lipschitz. Then for any $\mathcal{F} 
\in \ Lipflow(X,d)$ we have: 
$$f * \partial^{+} \mathcal{F} \ = \ \partial^{+}\left( f* \mathcal{F} \right) \ .$$
\label{lem3}
\end{lema}

This lemma finishes the proof of the theorem.  \quad $\blacksquare$
 
 We shall give now the proofs of the lemma \ref{lemacon}. The proofs of 
 lemmata \ref{lem2} and \ref{lem3} are almost straightforward, based on the unwinding of definitions.  
 
 \paragraph{Proof of Lemma \ref{lemacon}.}
 Let $\varepsilon > 0$ and $\mu \in (0,1)$. We are looking for $\varepsilon'  > 0$ and $\mu' \in (0,1)$ such that if $\displaystyle (x,y) \in V^{+}(\bar{\mathcal{F}}, \varepsilon', \mu')$ then $\displaystyle (y,x) \in V^{+}(\mathcal{F}, \varepsilon, \mu)$. 
 
 Choose $\displaystyle 4 \mu' < \min\left\{ \frac{1}{2}, \mu \right\}$. From condition (c) on the flow 
 $\mathcal{F}$ it follows that for any $t$ there is $\varepsilon(t) > 0$ such that for any $x,y \in X$, 
 $d(x,y) < \varepsilon(t)$ implies
 \begin{equation}
 \frac{1}{2 M(t)} d(\mathcal{F}(t,y), x) \ \leq \ d(y,\mathcal{F}(-t, x)) \ .
 \label{eq1lecon}
 \end{equation} 
 From conditions (a) and (c) we get that  for $\mid t \mid \leq \lambda_{0}$ we have 
 $M(t) \leq 2$. Let then $\varepsilon_{1} = \varepsilon(\lambda_{0})$.  The equation \eqref{eq1lecon} 
 implies that for any $x, y \in X$, $d(x,y) < \varepsilon_{1}$, and for any $\mid t \mid \leq \lambda_{0}$  
 we have: 
  \begin{equation}
 \frac{1}{4} d(\mathcal{F}(t,y), x) \ \leq \ d(y,\mathcal{F}(-t, x)) \ .
 \label{eq2lecon}
 \end{equation} 
 
 Let us choose $2 \sigma \leq \varepsilon_{1}$. From condition (b) it follows that there exists 
 $\varepsilon_{\sigma} > 0$ such that for any $\varepsilon" \in (0,\varepsilon_{\sigma})$ and any $x\in X$ 
 we have: 
 \begin{equation}
 d(x,\mathcal{F}(-\varepsilon", x) \ \leq \ \sigma \ . 
 \label{eq3lecon}
 \end{equation}
 Choose finally $\displaystyle \varepsilon' \leq \min \left\{ \varepsilon, \varepsilon_{\sigma}, \lambda_{0} \right\}$. Let us now take $\displaystyle (x,y) \in V^{+}(\bar{\mathcal{F}}, \varepsilon', \mu')$. This means that there is $\lambda' \in [0,\varepsilon']$ such that 
 \begin{equation}
 d(y, \mathcal{F}(-\lambda',x)) \ < \ \mu' d(x,y) \ . 
 \label{eq4lecon}
 \end{equation}
 A  consequence of inequalities \eqref{eq3lecon}, \eqref{eq4lecon} and the choices of $\varepsilon'$, $\mu'$  is that 
 $$d(x,y) \ < \ \frac{1}{1-\mu'} d(x, \mathcal{F}(-\lambda',x))  \ \leq \ \frac{\sigma}{1 -\mu'}  \ \leq \ 2 \sigma \leq \ \varepsilon_{1} \ .$$
 We can use then \eqref{eq2lecon} in \eqref{eq4lecon} to get: 
 $$\frac{1}{4} d(\mathcal{F}(\lambda, y), x) \ < \ \mu' d(x,y) \ .$$
 Use finally that $4 \mu' < \mu$ and $\lambda' \in [0,\varepsilon'] \subset [0,\varepsilon]$ and conclude: 
 there exists $\lambda \in [0,\varepsilon]$ such that: 
 $$d(\mathcal{F}(\lambda, y), x) \ < \ \mu d(x,y) \ .$$
 This means that $\displaystyle (y,x) \in V^{+}(\mathcal{F}, \varepsilon, \mu)$. \quad $\blacksquare$

\end{document}